\documentclass[12pt, leqno]{amsart}

\usepackage{amssymb,amsfonts,amscd,amsbsy}
\usepackage[mathscr]{eucal}

\theoremstyle{plain}
\newtheorem{thm}{Theorem}[section]
\newtheorem{lem}[thm]{Lemma}
\newtheorem{pro}[thm]{Proposition}
\newtheorem{cor}[thm]{Corollary}
\theoremstyle{definition}
\newtheorem{Def}[thm]{Definition}
\newtheorem{rem}[thm]{Remark}

\begin{document}

\title[Tame kernels and Further 4-rank densites]{Tame kernels and Further 4-rank densities}
\author[Robert Osburn and Brian Murray]{Robert Osburn and Brian Murray}

\address{Department of Mathematics $\&$ Statistics, McMaster University, 1280 Main Street West, Hamilton, Ontario L8S 4K1}

\address{Department of Mathematics, Louisiana State University,
Baton Rouge, LA 70803}

\email{osburnr@icarus.math.mcmaster.ca}
\email{murray@math.lsu.edu}


\begin{abstract}
There has been recent progress on computing the 4-rank of the tame kernel $K_2(\mathcal{O}_F)$ for $F$ a quadratic number field. For certain quadratic number fields, this progress has led to ``density results'' concerning the 4-rank of tame kernels. These results were first mentioned in \cite{CH} and proven in \cite{RO}. In this paper, we consider some additional quadratic number fields and obtain further density results of 4-ranks of tame kernels. Additionally, we give tables which might indicate densities in some generality. \end{abstract}

\maketitle
\section{Introduction}

We are interested in the structure of the 2-Sylow subgroup of
$K_2(\mathcal{O}_{F})$. As $K_{2}(\mathcal{O}_{F})$ is a finite abelian group, it is a
product of cyclic groups of prime power order.  We say the
$2^j$-rank, $j \geq 1$, of $K_{2}(\mathcal{O}_{F})$ is the
number of cyclic factors of $K_{2}(\mathcal{O}_{F})$ of order
divisible by $2^j$.  In \cite{Tate1}, the 2-rank of the tame
kernel is given by Tate's 2-rank formula.  In the case where F is
a quadratic number field, Browkin and Schinzel in \cite{BrS}
simplified the 2-rank formula.  In their formula, we can
determine the 2-rank by counting the number of elements in $\{\pm
1, \pm 2\}$ which are norms from the given quadratic field and
the number of odd primes which are ramified in the given
quadratic field. Now what about 4-ranks of $K_{2}(\mathcal{O}_F)$?

In \cite{CH}, Conner and Hurrelbrink characterize the 4-rank of $K_2(\mathcal{O})$ for certain quadratic number fields in terms of positive definite binary quadratic forms. This characterization led to a connection between densities of certain sets of primes and 4-rank values. Specifically,
the author in \cite{RO} considers the 4-rank of $K_2(\mathcal{O})$ for the quadratic number fields $\mathbb Q(\sqrt{pl})$, $\mathbb Q(\sqrt{2pl})$, $\mathbb Q(\sqrt{-pl})$, $\mathbb
Q(\sqrt{-2pl})$ for primes $p \equiv 7 \bmod 8$, $l \equiv 1 \bmod 8$ with $\Big(\frac{l}{p}\Big) =
1$.  In \cite{CH}, it was shown that for the fields $E = \mathbb Q(\sqrt{pl})$, $\mathbb
Q(\sqrt{2pl})$ and $F = \mathbb Q(\sqrt{-pl})$, $\mathbb Q(\sqrt{-2pl})$, $$ \mbox{4-rank}
\hspace{.05in} K_2(\mathcal{O}_E) = 1 \hspace{.05in} \mbox{or} \hspace{.05in} 2, $$ $$ \mbox{4-rank}
\hspace{.05in} K_2(\mathcal{O}_F) = 0 \hspace{.05in} \mbox{or} \hspace{.05in} 1. $$

The idea in \cite{RO} is to fix $p \equiv 7 \bmod 8$ and consider the set $$
\Omega = \{l \hspace{.05in}\mathrm{rational \hspace{.05in} prime}: l \equiv 1 \bmod 8 \hspace{.05in}
\mathrm{and} \hspace{.05in} \Big(\frac{l}{p}\Big) = \Big(\frac{p}{l}\Big) = 1 \}. $$

In \cite{RO}, the following was proved.

\begin{thm} \label{T:MQ0} For the fields $\mathbb Q(\sqrt{pl})$ and $\mathbb Q(\sqrt{2pl})$, 4-rank 1 and 2 each appear with natural density $\frac{1}{2}$ in
$\Omega$. For the fields $\mathbb Q(\sqrt{-pl})$ and $\mathbb Q(\sqrt{-2pl})$, 4-rank 0 and 1 each
appear with natural density $\frac{1}{2}$ in $\Omega$.
\end{thm}

In this paper, we consider the 4-rank of $K_2(\mathcal{O})$ for the quadratic number fields $\mathbb Q(\sqrt{pl})$, $\mathbb Q(\sqrt{-pl})$ for primes $p \equiv l \equiv 1 \bmod 8$ with $\Big(\frac{l}{p}\Big) = 1$ and $\mathbb Q(\sqrt{pl})$ for primes  $p \equiv l \equiv 1 \bmod 8$ with $\Big(\frac{l}{p}\Big) = -1$. We will see that for the primes  $p \equiv l \equiv 1 \bmod 8$ with $\Big(\frac{l}{p}\Big) = 1$,

\begin{center}
4-rank $K_2(\mathcal{O}_{\mathbb Q(\sqrt{pl})} ) = 1$ or $2$,
\end{center}
\begin{center}
4-rank $K_2(\mathcal{O}_{\mathbb Q(\sqrt{-pl})} ) = 1$ or $2$.
\end{center}

For the primes $p \equiv l \equiv 1 \bmod 8$ with $\Big(\frac{l}{p}\Big) = -1$, we will see
\begin{center}
4-rank $K_2(\mathcal{O}_{\mathbb Q(\sqrt{pl})} ) = 0$ or $1$.
\end{center}

Let us fix a prime $p \equiv 1 \bmod 8$ and consider the sets

$$
A = \{l \hspace{.05in}\mathrm{rational \hspace{.05in} prime}: l \equiv 1 \bmod 8 \hspace{.05in}
\mathrm{and} \hspace{.05in} \Big(\frac{l}{p}\Big) = 1 \}, $$

$$
B = \{l \hspace{.05in}\mathrm{rational \hspace{.05in} prime}: l \equiv 1 \bmod 8 \hspace{.05in}
\mathrm{and} \hspace{.05in} \Big(\frac{l}{p}\Big) = -1 \}. $$

The goal of this paper is to prove two theorems analogous to Theorem 1.1, namely:

\begin{thm} For the field $\mathbb Q(\sqrt{pl})$, 4-rank 1 and 2 appear with natural density $\frac{3}{4}$ and $\frac{1}{4}$ in $A$. For the field $\mathbb Q(\sqrt{-pl})$, 4-rank 1 and 2 each
appear with natural density $\frac{1}{2}$ in $A$.
\end{thm}

\begin{thm}
For the field $\mathbb Q(\sqrt{pl})$, 4-rank 0 and 1 each appear with natural density $\frac{1}{2}$ in $B$.
\end{thm}

Now for squarefree, odd integers d, consider the sets
$$X = \{d: d=pl \}$$
and
$$Y = \{d: d=-pl \}$$
for distinct primes $p$ and $l$.

We have computed the following: For $15 \le d < 10^6$, there are $168331$ d's in $X$. Among them, there are $35787$ d's (21.26$\%$) yielding 4-rank 0, $128468$ d's (76.32$\%$) yielding 4-rank 1, and $4076$ d's (2.42$\%$) yielding 4-rank 2.

For $-10^6 < d \le -15$, there are $168330$ d's in $Y$. Among them, there are
$104056$ d's (61.82$\%$) yielding 4-rank 0, $63054$ d's (37.46$\%$) yielding 4-rank 1, and $1220$ d's (.72$\%$) yielding 4-rank 2. As a consequence of Theorems 1.2, 1.3 and Tables I and II in \cite{Qin1} and \cite{Qin2}, we obtain:

\begin{cor}
For the fields $\mathbb Q(\sqrt{pl})$, 4-rank 0, 1, and 2 appear with natural density $\frac{13}{64}$, $\frac{97}{128}$, $\frac{5}{128}$ respectively in $X$.
\end{cor}

\begin{cor} For the fields $\mathbb Q(\sqrt{-pl})$, 4-rank 0, 1, and 2 appear with natural density $\frac{37}{64}$, $\frac{13}{32}$, and $\frac{1}{64}$ respectively in $Y$.
\end{cor}

\section{Preliminaries}

Let $\mathcal{D}$ be a Galois extension of $\mathbb Q$, and $G = Gal(\mathcal{D}/ \mathbb Q)$.  Let
$Z(G)$ denote the center of G and $\mathcal{D}^{Z(G)}$ denote the fixed field of $Z(G)$.  Let $p$ be a
rational prime which is unramified in $\mathcal{D}$ and $\beta$ be a prime of $\mathcal{D}$ containing
$p$.  Let $\Big(\frac{\mathcal{D}/\mathbb Q}{p}\Big)$ denote the Artin symbol of $p$ and $\{g\}$ the
conjugacy class containing one element $g \in G$. In Sections 5 and 6 we use the following elementary lemma from \cite{RO}.

\begin{lem} \label{L:center} $\Big(\frac{\mathcal{D}/\mathbb Q}{p}\Big) = \{g\}$ for some $g \in Z(G)$ if and only if p splits completely
in $\mathcal{D}^{Z(G)}$.
\end{lem}

Thus if we can show that rational primes split completely in the fixed field of the center of a certain
Galois group G, then we know the associated Artin symbol is a conjugacy class containing one element. Note that determining the order of $Z(G)$ gives us the number of possible
choices for the Artin symbol. The order of $Z(G)$ can be computed using the following setup.

Let $G_1$ and $G_2$ be finite groups and $A$ a finite abelian group. Suppose $r_1:G_1 \to A$ and
$r_2:G_2 \to A$ are two epimorphisms and $\mathcal{G} \subset G_1 \times G_2$ is the set $\{(g_1,
g_2)\in G_1 \times G_2 : r_1(g_1)=r_2(g_2) \}$.  Since A is abelian, there is an epimorphism $r:G_1
\times G_2 \to A$ given by $r(g_1, g_2)=r_1(g_1){r_2(g_2)}^{-1}$.  Thus $\mathcal{G} = ker(r) \subset
G_1 \times G_2$. One can check that $Z(\mathcal{G}) = \mathcal{G} \cap Z(G_1 \times G_2)$. From \cite{RO}, we provide:

\begin{lem} \label{L:easy}
$Z(\mathcal{G}) = Z(G_1) \times Z(G_2) \iff r_1\Big |_{Z(G_1)}$ and $r_2\Big |_{Z(G_2)}$ are both trivial.
\end{lem}

We will use the following definition throughout this paper.

\begin{Def} \label{D:sat} For primes $p \equiv l \equiv 1 \bmod 8$ with $\Big(\frac{l}{p}\Big) = 1$, $\mathcal{K} = \mathbb Q$$(\sqrt{2p})$, and $h^{+}(\mathcal{K})$ the narrow class number of $\mathcal{K}$, we say:

{\em l satisfies $\langle 1, 32\rangle$} if and only if $l = x^2 + 32y^2$ for some $x,y \in \mathbb Z$

{\em l satisfies $\langle p, -2\rangle$} if and only if $l^{\frac{h^{+}(\mathcal{K})}{4}} = pn^{2} - 2m^{2}$
for some $n,m \in \mathbb Z$ with $m \not \equiv 0\bmod l$

{\em l satisfies $\langle 1, -2p\rangle$} if and only if $l^{\frac{h^{+}(\mathcal{K})}{4}} = n^{2} - 2pm^{2}$ for some
$n,m \in \mathbb Z$ with $m \not \equiv 0\bmod l$.
\end{Def}

\section{First Extension}
Consider the fixed prime $p \equiv 1\bmod 8$.  Note $p$ splits completely in $\mathcal{L} = \mathbb
Q(\sqrt{2})$ over $\mathbb Q$ and so $$ p\mathcal{O}_{\mathcal{L}} = \frak B \frak B^{'} $$ for some
primes $\frak B \not = \frak B^{'}$ in $\mathcal{L}$.  The field $\mathcal{L}$ has narrow class number
$h^{+}(\mathcal{L}) = 1$ as $h(\mathcal{L}) = 1$ and $N_{\mathcal{L}/\mathbb Q}(\epsilon) = -1$ where
$\epsilon = 1 + \sqrt{2}$ is a fundamental unit of $\mathbb Q(\sqrt{2})$.
Similar to Lemma 2.1 in \cite{CH},

\begin{lem} \label{L:pi} The prime $\frak B$ which occurs in the decomposition of $p\mathcal{O}_{\mathcal{L}}$ has a generator $\pi = a + b\sqrt{2} \in \mathcal{O}_{\mathcal{L}}$, unique up to a sign and to multiplication by the square of a unit in $\mathcal{O}_{\mathcal{L}}^{*}$ for which $N_{\mathcal{L}/\mathbb Q}(\pi) = a^2 - 2b^2 = p$.
\end{lem}

The degree 4 extension $\mathbb Q(\sqrt{2}, \sqrt{\pi})$ over $\mathbb Q$ has normal closure

\noindent{$\mathbb Q(\sqrt{2}, \sqrt{\pi}, \sqrt{p})$} as $N_{\mathcal{L}/\mathbb Q}(\pi) = p$.  Set
$$ N = \mathbb Q(\sqrt{2}, \sqrt{\pi}, \sqrt{p}). $$ Then $N$ is Galois over $\mathbb Q$ and $[N: \mathbb Q] = 8$.  By Corollary 24.5 in \cite{CNP}, $4$ divides the narrow class number of $\mathbb Q(\sqrt{2p})$. Moreover $N$ over $\mathbb Q(\sqrt{2p})$ is unramified at all finite primes. Similar to Lemma 2.3
in \cite{CH}, $N$ is the unique unramified cyclic degree 4 extension over $\mathbb Q(\sqrt{2p})$.

Consider the rational primes $l \equiv 1 \bmod 8$ for which $\Big(\frac{l}{p}\Big) = 1$. These primes split completely in $\mathbb Q(\sqrt{2}, \sqrt{p})$ over $\mathbb Q$. We characterize such primes $l$ that split completely in $N$ over $\mathbb Q$. As $N$ is the unique unramified cyclic degree 4 extension of $\mathbb Q(\sqrt{2p})$, mimicing Lemma 3.3 in \cite{CH} yields

\begin{lem} Let $l \equiv 1 \bmod 8$ be a prime such that $\Big(\frac{l}{p}\Big) = 1$. Then:

$l$ splits completely in $N$ if and only if $l$ satisfies $\langle 1, -2p\rangle$.
\end{lem}

Similar to Lemma 3.4 in \cite{CH}, with 2 (respectively, $\frak{D}$, the unique dyadic prime in $\mathcal{O}_{\mathbb Q(\sqrt{2p})}$ ) replaced by $p$ (respectively $\frak{p}$, the prime over $p$ whose class is the unique element of order 2 in the narrow ideal class group of $\mathbb Q(\sqrt{2p})$ ), we obtain

\begin{lem} Let $l \equiv 1 \bmod 8$ be a prime such that $\Big(\frac{l}{p}\Big) = 1$. Then:

$l$ does not split completely in $N$ if and only if $l$ satisfies $\langle p, -2\rangle$.
\end{lem}

We now relate the characterizations of Lemmas 3.2 and 3.3 to the quadratic symbol $\Big(\frac{\pi}{l}\Big)$. From Lemma 3.1, we have a presentation $\pi = a + b\sqrt{2} \in \mathcal{O}_{\mathcal{L}}$ with $N_{\mathcal{L}/\mathbb Q}(\pi) = p$. Let $\frak P$ be a prime above $l$ in $\mathcal{O}_{\mathcal{L}}$. As $l$ splits in $\mathcal{L}$ over $\mathbb Q$, then the residue field $\mathcal{O}_{\mathcal{L}} / \frak P$ is isomorphic to $\mathbb Z /{l}\mathbb Z = {\mathbb F}_{l}$, the field with $l$ elements. As 2 is a square modulo $l$, we have $2 \equiv \alpha^{2} \bmod l$ for some $\alpha \in {{\mathbb F}_{l}}^{*}$.
Thus we can identify $\pi = a + b\sqrt{2} \in \mathcal{O}_{\mathcal{L}}$
with $a + b\alpha \in {\mathbb F}_{l}$. When we write the symbol $\Big(\frac{\pi}{l}\Big)$, it is understood that we mean $\Big(\frac{a + b \alpha}{l}\Big)$. From the discussion in Section 3 of \cite{CH}, the symbol
$\Big(\frac{\pi}{l}\Big)$ is well defined and $l$ splits completely in $N$ over $\mathbb Q$ if and only if $\Big(\frac{\pi}{l}\Big) = 1$. Combining this discussion with Lemmas 3.2 and 3.3, we have:

\begin{pro}  Let $l \equiv 1 \bmod 8$ be a prime with $\Big(\frac{l}{p}\Big) = 1$. Then:
\begin{center}
$l$ satisfies $\langle 1, -2p\rangle$ $\iff$ $\Big(\frac{\pi}{l}\Big) = 1$,

$l$ satisfies $\langle p, -2\rangle$ $\iff$ $\Big(\frac{\pi}{l}\Big) = -1$.
\end{center}
\end{pro}

\section{Matrices and Symbols}
Hurrelbrink and Kolster \cite{HK98} generalize Qin's approach in \cite{Qin1}, \cite{Qin2} and obtain 4-rank results
by computing $\Bbb F_2$-ranks of certain matrices of local
Hilbert symbols. Let us be more specific. Let $F=\mathbb Q(\sqrt{d})$, $d\neq 0,1$, squarefree.  Let $p_1, p_2, \dots , p_t$ denote the odd primes dividing $d$. Recall 2 is a norm from F $\iff$ all $p_{i}$'s are $\equiv \pm 1 \bmod 8.$ If so, then d is a norm from $\mathbb Q(\sqrt{2})$, thus
$$
d=u^2-2w^2
$$
for $u,w \in \mathbb Z$. Now consider two matrices:

If $d<0$,

${M^{'}}_{F/\mathbb Q}=
\left( \begin{matrix} (-d,p_1)_2 & (-d,p_1)_{p_1} \hspace{.05in} \dots & (-d,p_1)_{p_t} \\
(-d,p_2)_2 & (-d,p_2)_{p_1} \hspace{.05in} \dots & (-d,p_2)_{p_t} \\
\vdots & \vdots & \vdots \\
(-d,p_{t-1})_2 & (-d,p_{t-1})_{p_1} \hspace{.05in} \dots & (-d,p_{t-1})_{p_t} \\
(-d,v)_2 & (-d,v)_{p_1} \hspace{.05in} \dots & (-d,v)_{p_t} \\
(-d,-1)_2 & (-d,-1)_{p_1} \hspace{.05in} \dots & (-d,-1)_{p_t} \\
\end{matrix} \right)
$

If $d>0$,

$M_{F/\mathbb Q}=
\left( \begin{matrix} (-d,p_1)_2 & (-d,p_1)_{p_1} \hspace{.1in} \dots & (-d,p_1)_{p_t} \\
(-d,p_2)_2 & (-d,p_2)_{p_1} \hspace{.1in} \dots & (-d,p_2)_{p_t} \\
\vdots & \vdots & \vdots \\
(-d,p_{t-1})_2 & (-d,p_{t-1})_{p_1} \hspace{.1in} \dots & (-d,p_{t-1})_{p_t} \\
(-d,v)_2 & (-d,v)_{p_1} \hspace{.1in} \dots & (-d,v)_{p_t} \\
(d,-1)_2 & (d,-1)_{p_1} \hspace{.1in} \dots & (d,-1)_{p_t} \\
\end{matrix} \right)$

If 2 is not a norm from F, set $v=2$. Otherwise, set $v=u+w$.
Replacing the $1$'s by $0$'s and the $-1$'s by $1$'s, we calculate
the matrix rank over $\mathbb F_2$. Why look at these matrices ?
From \cite{HK98},

\begin{lem} Let $F=\mathbb Q(\sqrt{d})$, $d \neq 0,1$, squarefree. Then \\
(i)If $d<0$, then 4-rank $K_2(\mathcal{O}_F) = t$ $-$ rk $({M^{'}}_{F/\mathbb Q})$ \\
(ii) If $d>0$, then 4-rank $K_2(\mathcal{O}_F) = t$ $-$ rk $(M_{F/\mathbb Q})$ + $a^{'} - a$
 \\

where \\

$a=
\left \{ \begin{array}{l}
0 \quad \mbox{if 2 is a norm from F} \\
1 \quad \mbox{otherwise}
\end{array}
\right. $ \\

and \\

$a^{'}=
\left \{ \begin{array}{l}
0 \quad \mbox{if both -1 and 2 are norms from F} \\
1 \quad \mbox{if exactly one of -1 or 2 is a norm from F} \\
2 \quad \mbox{if none of -1 or 2 are norms from F.} \\
\end{array}
\right. $

\end{lem}

Recall that our cases are:
\begin{itemize}
\item $\mathbb Q(\sqrt{pl})$, $\mathbb Q(\sqrt{-pl})$ where  $p \equiv l \equiv 1 \bmod 8$ with $\Big(\frac{l}{p}\Big) = 1$, \\
\item $\mathbb Q(\sqrt{pl})$ for  $p \equiv l \equiv 1 \bmod 8$ with $\Big(\frac{l}{p}\Big) = -1$.
\end{itemize}

In both cases 2 is a norm from $\mathbb Q(\sqrt{pl})$ and $\mathbb Q(\sqrt{-pl})$. Before we view the matrices for our cases, we characterize the symbol $(-d,v)_2$ for $d=pl, -pl$ (see Lemmas 5.3 amd 5.15 in \cite{HK98}).

\begin{itemize}
\item $(-pl,v)_2 = 1 \iff$ both $p$, $l$ satisfy $\langle 1, 32\rangle$ or neither $p$, $l$ satisfy $\langle 1, 32\rangle$, \\
\item $(pl,v)_2 = 1$.
\end{itemize}

Also, $v$ is an $l$-adic unit and hence
\begin{center}
$(-pl,v)_{l} = (l,v)_{l} = \Big(\frac{v}{l}\Big)$.
\end{center}

Similarly, $(-pl,v)_p = \Big(\frac{v}{p}\Big)$. In the entries of the matrices below, we write $(-pl,v)_2$, $\Big(\frac{v}{l}\Big)$, and $\Big(\frac{v}{p}\Big)$ remembering to first evalute the symbols, make the substitutions 1 for 0 and -1 for 1, and then calculate the matrix rank over $\mathbb F_2$. Now what are the matrices in our situations?

\begin{itemize}
\item For $p \equiv l \equiv 1 \bmod 8$ with $\Big(\frac{l}{p}\Big) = 1$, we have:

$M_{\mathbb Q(\sqrt{pl})/\mathbb Q}=
\left( \begin{matrix} 0 & 0 & 0 \\
(-pl,v)_2 & \Big(\frac{v}{p}\Big) & \Big(\frac{v}{l}\Big) \\
0 & 0 & 0 \\
\end{matrix} \right)$,

${M^{'}}_{\mathbb Q(\sqrt{-pl})/\mathbb Q}=
\left( \begin{matrix} 0 & 0 & 0 \\
0 & \Big(\frac{v}{p}\Big) & \Big(\frac{v}{l}\Big) \\
0 & 0 & 0 \\
\end{matrix} \right)$.

\item For $p \equiv l \equiv 1 \bmod 8$ with $\Big(\frac{l}{p}\Big) = -1$, we have:

$M_{\mathbb Q(\sqrt{pl})/\mathbb Q}=
\left( \begin{matrix} 0 & 1 & 1 \\
(-pl,v)_2 & \Big(\frac{v}{p}\Big) & \Big(\frac{v}{l}\Big) \\
0 & 0 & 0 \\
\end{matrix} \right)$.

\end{itemize}

\begin{rem} For $p \equiv l \equiv 1 \bmod 8$ with $\Big(\frac{l}{p}\Big) = 1$, we have:
\begin{itemize}
\item 4-rank $K_2(\mathcal{O}_{\mathbb Q(\sqrt{pl})}) = 1
\iff$ rank $M_{\mathbb Q(\sqrt{pl})/\mathbb Q} = 1
\iff (-pl,v)_2 = 1$, $\Big(\frac{v}{l}\Big) = -1$ or $(-pl,v)_2 = -1
\iff$ both $p$, $l$ satisfy $\langle 1, 32\rangle$, $\Big(\frac{v}{l}\Big) = -1$
or neither $p$, $l$ satisfy $\langle 1, 32\rangle$, $\Big(\frac{v}{l}\Big) = -1$, or exactly one of $p$, $l$ satisifes $\langle 1, 32\rangle$. \\

\item 4-rank $K_2(\mathcal{O}_{\mathbb Q(\sqrt{pl})}) = 2
\iff$ rank $M_{\mathbb Q(\sqrt{pl})/\mathbb Q} = 0
\iff (-pl,v)_2 = 1$, $\Big(\frac{v}{l}\Big) = 1
\iff$ both $p$, $l$ satisfy $\langle 1, 32\rangle$, $\Big(\frac{v}{l}\Big) = 1$
or neither $p$, $l$ satisfy $\langle 1, 32\rangle$, $\Big(\frac{v}{l}\Big) = 1$.\\

\item 4-rank $K_2(\mathcal{O}_{\mathbb Q(\sqrt{-pl})}) = 1
\iff$ rank ${M^{'}}_{\mathbb Q(\sqrt{-pl})/\mathbb Q} = 1
\iff \Big(\frac{v}{l}\Big) = -1$.

\item  4-rank $K_2(\mathcal{O}_{\mathbb Q(\sqrt{-pl})}) = 2
\iff$ rank ${M^{'}}_{\mathbb Q(\sqrt{-pl})/\mathbb Q} = 0
\iff \Big(\frac{v}{l}\Big) = 1$.

\end{itemize}
\end{rem}

\begin{rem} For $p \equiv l \equiv 1 \bmod 8$ with $\Big(\frac{l}{p}\Big) = -1$:
\begin{itemize}
\item 4-rank $K_2(\mathcal{O}_{\mathbb Q(\sqrt{pl})}) = 1
\iff$ rank $M_{\mathbb Q(\sqrt{pl})/\mathbb Q} = 1
\iff (-pl,v)_2 = 1
\iff$ both $p$, $l$ satisfy $\langle 1, 32\rangle$ or neither $p$, $l$ satisfy $\langle 1, 32\rangle$. \\

\item  4-rank $K_2(\mathcal{O}_{\mathbb Q(\sqrt{pl})}) = 0
\iff$ rank $M_{\mathbb Q(\sqrt{pl})/\mathbb Q} = 2
\iff (-pl,v)_2 = -1
\iff$ exactly one of $p$, $l$ satisfies $\langle 1, 32\rangle$.
\end{itemize}
\end{rem}

We can now prove Theorem 1.3.

\begin{proof}
Consider the sets
\begin{center}
$\mathcal{A}_1 = \{l$ prime: $l \equiv 1 \bmod 8$ and $l$ satisfies $\langle 1, 32\rangle$ \},
\end{center}
\begin{center}
$\mathcal{A}_2 = \{l$ prime: $l \equiv 1 \bmod 8$ and $l$ does not satisfy $\langle 1, 32\rangle$ \}.
\end{center}
By the discussion before Corollary 24.2 in \cite{CNP}, $\mathcal{A}_1$ and $\mathcal{A}_2$ each have density $\frac{1}{2}$ in the set of all primes $l \equiv 1 \bmod 8$. By Dirichlet's Theorem on primes in arithmetic progessions, $\mathcal{A}_1$ and $\mathcal{A}_2$ each have density $\frac{1}{8}$ in the set of all primes $l$. Note that for primes $p \equiv 1 \bmod 8$, the sets
\begin{center}
$\mathcal{B}_1 = \{l$ prime: $l \equiv 1 \bmod 8$, $\Big(\frac{l}{p}\Big) = -1$,  and $l$ satisfies $\langle 1, 32\rangle$ \},
\end{center}
\begin{center}
$\mathcal{B}_2 = \{l$ prime: $l \equiv 1 \bmod 8$, $\Big(\frac{l}{p}\Big) = -1$, and $l$ does not satisfy $\langle 1, 32\rangle$ \}.
\end{center}

\noindent each have density $\frac{1}{2}$ in $\mathcal{A}_1$ and $\mathcal{A}_2$ respectively. Thus $\mathcal{B}_{1}$ and $\mathcal{B}_{2}$ have densities $\frac{1}{16}$ in the set of all primes $l$. If $p$ satisfies $\langle 1, 32\rangle$, then by Remark 4.3:
\begin{center}
$\mathcal{B}_1 = \{l$ prime: $l \equiv 1 \bmod 8$, $\Big(\frac{l}{p}\Big) = -1$,  and 4-rank $K_2(\mathcal{O}_{\mathbb Q(\sqrt{pl})}) = 1$ \},
\end{center}
\begin{center}
$\mathcal{B}_2 = \{l$ prime: $l \equiv 1 \bmod 8$, $\Big(\frac{l}{p}\Big) = -1$, and 4-rank $K_2(\mathcal{O}_{\mathbb Q(\sqrt{pl})}) = 0$\}.
\end{center}

For each $\mathcal{B}_{i}$, $i=1,2$, we have:

$\left \{ \begin{array}{l}
\mbox{Density of $\mathcal{B}_{i}$ in the} \\
\mbox{set of all primes $l$}
\end{array}
\right \}$ $=$ $\left \{ \begin{array}{l}
\mbox{Density of} \\
\mbox{$\mathcal{B}_{i}$ in $B$}
\end{array}
\right \} $ $\cdot$ $\left \{ \begin{array}{l}
\mbox{Density of $B$ in the} \\
\mbox{set of all primes $l$}
\end{array}
\right \}$

\noindent where $B$ has density $\frac{1}{8}$ in the set of all primes $l$. Thus 4-rank 0 and 4-rank 1 each appear with natural density $\frac{1}{2}$ in $B$. A similar argument works if $p$ does not satisfy $\langle 1, 32\rangle$.

\end{proof}

For the primes $p \equiv l \equiv 1 \bmod 8$ with $\Big(\frac{l}{p}\Big) = 1$, let us relate the Legendre symbol $\Big(\frac{v}{l}\Big)$ to the quadratic symbol $\Big(\frac{\pi}{l}\Big)$. For primes $l \equiv 1 \bmod 8$, the quadratic symbol $\Big(\frac{1+\sqrt{2}}{l}\Big)$ is well defined and satisfies, see \cite{BC},
\begin{center}
$\Big(\frac{1+\sqrt{2}}{l}\Big) = 1 \iff l$ satisfies $\langle 1, 32\rangle$.
\end{center}

\begin{pro} Let $d=\pm pl$ be as above, $d=u^2-2w^2$ with $u$,$w \in \mathbb Z$.
Then:
\begin{center}
$\Big(\frac{v}{l}\Big) = \Big(\frac{\pi}{l}\Big) \Big(\frac{1+\sqrt{2}}{l}\Big)$ if $d=pl$
\end{center}
\begin{center}
$\Big(\frac{v}{l}\Big) = \Big(\frac{\pi}{l}\Big)$   if $d=-pl$.
\end{center}
\end{pro}

\begin{proof}
From the proof of Proposition 4.6 in \cite{CH}, we use the identity
\begin{center}
$\Big(\frac{v}{l}\Big) = \Big(\frac{\gamma + \delta\sqrt{2}}{l}\Big) \Big(\frac{1+\sqrt{2}}{l}\Big)$
\end{center}

\noindent where $\frac{d}{l} =  N_{\mathcal{L}/\mathbb Q}(\gamma + \delta\sqrt{2})$ for $\gamma$, $\delta \in \mathbb Z$. For $d=pl$, we have $\frac{d}{l} = p =  N_{\mathcal{L}/\mathbb Q}(\pi)$ and thus $\gamma + \delta\sqrt{2} = \pi$, up to squares. For $d=-pl$, we have  $\frac{d}{l} = -p =  -N_{\mathcal{L}/\mathbb Q}(\pi)$ and so $\gamma + \delta\sqrt{2} = (1+\sqrt{2})\pi$, up to squares.

\end{proof}

In view of Proposition 3.4, Remark 4.2, and Proposition 4.4, we can determine the 4-rank of the tame kernel in terms of quadratic forms.

\begin{pro} For $p \equiv l \equiv 1 \bmod 8$ with $\Big(\frac{l}{p}\Big)=1$:
\begin{itemize}
\item 4-rank $K_2(\mathcal{O}_{\mathbb Q(\sqrt{pl})}) = 1 \iff$ both $p$, $l$ satisfy $\langle 1, 32\rangle$, $l$ satisfies $\langle p, -2\rangle$ or neither $p$, $l$ satisfy $\langle 1, 32\rangle$, $l$ satisfies $\langle p, -2\rangle$ or exactly one of $p$, $l$ satisfies $\langle 1, 32\rangle$

\item 4-rank $K_2(\mathcal{O}_{\mathbb Q(\sqrt{pl})}) = 2 \iff$ both $p$, $l$ satisfy $\langle 1, 32\rangle$, $l$ satisfies $\langle 1, -2p\rangle$ or neither $p$, $l$ satisfy $\langle 1, 32\rangle$, $l$ satisfies $\langle 1, -2p\rangle$

\item 4-rank $K_2(\mathcal{O}_{\mathbb Q(\sqrt{-pl})}) = 1 \iff l$ satisfies $\langle p, -2\rangle$

\item  4-rank $K_2(\mathcal{O}_{\mathbb Q(\sqrt{-pl})}) = 2 \iff l$ satisfies $\langle 1, -2p\rangle$.
\end{itemize}
\end{pro}

It should be noted that Qin Yue has obtained characterizations of
4-rank values, similar to Proposition 4.5, by additionally
assuming that the fundamental unit of $\mathbb Q(\sqrt{2p})$, $p
\equiv 1 \bmod 8$, has norm $-1$, see \cite{Yue}.

\section{Two Artin Symbols}
\subsection{First Artin symbol}
Consider $\mathbb Q(\sqrt{2})$ over $\mathbb Q$.  Let $\epsilon = 1 + \sqrt{2} \in (\mathbb
Z[\sqrt{2}])^{*}$. Then $\epsilon$ is a fundamental unit of $\mathbb Q(\sqrt{2})$ which has norm $-1$.
The degree 4 extension $\mathbb Q(\sqrt{2}, \sqrt{\epsilon})$ over $\mathbb Q$ has normal closure
$\mathbb Q(\sqrt{2}, \sqrt{\epsilon}, \sqrt{-1})$.  Set $$ N_1 = \mathbb Q(\sqrt{2}, \sqrt{\epsilon},
\sqrt{-1}).$$ Note that $Gal(N_1/\mathbb Q)$ is the dihedral group of order 8  and $Z(Gal(N_1/\mathbb Q)) = Gal(N_1/\mathbb Q(\sqrt{2}, \sqrt{-1}))$ (see \cite{RO}, Section 3.2).

Only the prime 2 ramifies in $\mathbb Q(\sqrt{2})$, $\mathbb Q(\sqrt{-1})$, $\mathbb
Q(\sqrt{\epsilon})$, and so only the prime 2 ramifies in the compositum $N_1$ over $\mathbb Q$. Now as
$l \in A$ is unramified in $N_1$ over $\mathbb Q$, the Artin symbol $\Big(\frac{N_1/\mathbb Q}{\beta}\Big)$ is defined for primes $\beta$ of $\mathcal{O}_{N_1}$ containing $l$.  Let
$\Big(\frac{N_1/\mathbb Q}{l}\Big)$ denote the conjugacy class of $\Big(\frac{N_1/\mathbb Q}{\beta}\Big)$ in $Gal(N_1/\mathbb Q)$.  The primes ${l \in A}$ split completely in $\mathbb Q(\sqrt{2},
\sqrt{-1})$ and ${N_1}^{Z(Gal(N_1/\mathbb Q))}  = \mathbb Q(\sqrt{2}, \sqrt{-1}).$ Thus by Lemma
\ref{L:center}, we have that $\Big(\frac{N_1/\mathbb Q}{l}\Big) = \{g\}$ for some $g \in Z(Gal(N_1/\mathbb Q))$. As $Z(Gal(N_1/\mathbb Q)$) has order 2, there are two possible choices for $\Big(\frac{N_1/\mathbb Q}{l}\Big)$.  Combining this statement with Addendum (3.7) from \cite{CH}, we have

\begin{rem} \label{R:N1}

\begin{eqnarray}
\Big(\frac{N_1/\mathbb Q}{l}\Big) = \{id\}
& \iff & l \hspace{.05in} \mbox{splits completely in} \hspace{.05in} N_1  \nonumber \\
& \iff & l \hspace{.05in} \mbox{satisfies} \hspace{.05in} \langle 1, 32\rangle. \nonumber
\end{eqnarray}

\end{rem}

\subsection{Second Artin symbol}
In section 3, we considered
\begin{center}
$N= \mathbb Q(\sqrt{2}, \sqrt{\pi}, \sqrt{p})$,
\end{center}
the unique unramified cyclic degree 4 extension over $\mathbb Q(\sqrt{2p})$. Similar to the extension $N_1$, we have $Gal(N/\mathbb Q)$ is the dihedral group of order 8 and $Z(Gal(N/\mathbb Q)) = Gal(N/\mathbb Q(\sqrt{2}, \sqrt{p}))$.

\begin{pro} \label{P:ramN2} If $l \in A$, then l is unramified in $N$ over $\mathbb Q$.
\end{pro}

\begin{proof} Since $p \equiv 1\bmod 8$, the discriminant of $\mathbb Q(\sqrt{2p})$ is $8p$.  For $l \in A$, we have $\Big(\frac{2p}{l}\Big) = 1$ and so $l$ is unramified in $\mathbb Q(\sqrt{2p})$.  We conclude that $l$ is unramified in $N$ over $\mathbb Q$.
\end{proof}

As $l \in A$ is unramified in $N$ over $\mathbb Q$, the Artin symbol $\Big(\frac{N/\mathbb Q}{\beta}\Big)$ is defined for primes $\beta$ of $\mathcal{O}_{N}$ containing $l$.  Let
$\Big(\frac{N/\mathbb Q}{l}\Big)$ denote the conjugacy class of $\Big(\frac{N/\mathbb Q}{\beta}\Big)$ in $Gal(N/\mathbb Q)$.  The primes $l \in A$ split completely in $\mathbb Q(\sqrt{2},
\sqrt{p})$ and ${N}^{Z(Gal(N/\mathbb Q))} = \mathbb Q(\sqrt{2}, \sqrt{p}).$ By Lemma
\ref{L:center}, we have that $\Big(\frac{N/\mathbb Q}{l}\Big) = \{h\}$ for some $h \in Z(Gal(N/\mathbb Q))$.  As Z(Gal$(N/\mathbb Q))$ has order 2, there are two
possible choices for $\Big(\frac{N/\mathbb Q}{l}\Big)$.  Combining this statement and Lemmas 3.2 and 3.3, we have

\begin{rem} \label{R:N2}

\begin{eqnarray}
\Big(\frac{N/\mathbb Q}{l}\Big) = \{id\}
& \iff & l \hspace{.05in} \mbox{splits completely in} \hspace{.05in} N \nonumber \\
& \iff & l \hspace{.05in} \mbox{satisfies} \hspace{.05in} \langle 1, -2p\rangle. \nonumber
\end{eqnarray}

\begin{eqnarray}
\Big(\frac{N/\mathbb Q}{l}\Big) \not = \{id\}
& \iff & l \hspace{.05in} \mbox{does not split completely in} \hspace{.05in} N \nonumber \\
& \iff & l \hspace{.05in} \mbox{satisfies} \hspace{.05in} \langle p, -2\rangle. \nonumber
\end{eqnarray}

\end{rem}

\section{A composite and proof of Theorem 1.2}
In this section we consider the composite field $N_1N$. Set $$ \frak N = N_1N. $$ Note that $[\frak N:\mathbb Q] = 32$. As $N_1$ and $N$ are normal extensions of $\mathbb Q$, $\frak N$ is a normal extension of $\mathbb Q$.

For $l \in A$, $l$ is unramified in $\frak N$ as it is unramified in $N_1$ and $N$.  The Artin symbol $\Big(\frac{\frak N/\mathbb Q}{\beta}\Big)$ is now defined for some
prime $\beta$ of $\mathcal{O}_{\frak N}$ containing $l$.  Let $\Big(\frac{\frak N/\mathbb Q}{l}\Big)$ denote the
conjugacy class of $\Big(\frac{\frak N/\mathbb Q}{\beta}\Big)$ in $Gal(\frak N/\mathbb Q)$.  Letting $M =
\mathbb Q(\sqrt{2}, \sqrt{-1}, \sqrt{p}) \subset \frak N$, we prove

\begin{lem} \label{L:N_1Ncenter} $Z(Gal({\frak N}/\mathbb Q)) = Gal({\frak N}/M)$ is elementary abelian of order 4.
\end{lem}

\begin{proof} For $\sigma \in Gal({\frak N}/M)$, $\sigma$ can only change the sign of $\sqrt{\epsilon}$ and $\sqrt{\pi}$ as $\epsilon \in M$. Since $\frak N =M(\sqrt{\epsilon}, \sqrt{\pi})$, $Gal({\frak N}/M)$ is elementary abelian of order 4. Now consider the restrictions $r_1:G_1 \to Gal(\mathbb Q(\sqrt{2})/\mathbb Q)$ and $r_2:G_2 \to
Gal(\mathbb Q(\sqrt{2})/\mathbb Q)$ where $G_1 = Gal(N_1/\mathbb Q)$ and $G_2 = Gal(N/\mathbb Q)$.
Clearly $r_1\Big |_{Z(G_1)}$ and $r_1\Big |_{Z(G_2)}$ are both trivial.  Then by Lemma
\ref{L:easy}, $Z(\mathcal{G})$ is elementary abelian of order 4 where $\mathcal{G} =
Gal({\frak N}/\mathbb Q)$. Thus
$Z(Gal({\frak N}/\mathbb Q)) = Gal({\frak N}/M)$.

\end{proof}

Now for $l \in A$, $l$ splits completely in $\mathbb Q(\sqrt{-1})$ and $\mathbb Q(\sqrt{2},\sqrt{p})$ and so splits completely in the composite field $M = \mathbb Q(\sqrt{2}, \sqrt{-1}, \sqrt{p})$.  From Lemma \ref{L:N_1Ncenter}, ${\frak N}^{Z(Gal({\frak N}/\mathbb Q))} = \mathbb Q(\sqrt{2}, \sqrt{-1},
\sqrt{p}).$ So by Lemma \ref{L:center}, we have
\begin{center}
$\Big(\frac{{\frak N}/\mathbb Q}{l}\Big) = \{k\}$ for some $k \in Gal({\frak N}/\mathbb Q)$.
\end{center}

As $Z(Gal({\frak N}/\mathbb Q))$ has order 4, there
are four possible choices for $\Big(\frac{{\frak N}/\mathbb Q}{l}\Big)$.  Using Remarks 5.1 and 5.3, we now make the following one to one correspondences.

\begin{rem} \label{R:four}
(i) $\Big(\frac{{\frak N}/\mathbb Q}{l}\Big) = \{id\} \iff$ $l$ splits completely in $\frak N$  $\iff$ $\left \{
\begin{array}{l}
l \hspace{.05in} \mbox{splits completely in} \\
\mbox{$N_1$ and $N$}
\end{array}
\right \}$ $\iff$ $\left \{ \begin{array}{l}
l \hspace{.05in} \mbox{satisfies} \hspace{.05in} \langle 1, 32\rangle \\
l \hspace{.05in} \mbox{satisfies} \hspace{.05in} \langle 1, -2p\rangle \\
\end{array}
\right \} $.

(ii)  $\Big(\frac{{\frak N}/\mathbb Q}{l}\Big) \not = \{id\} \iff$ $l$ does not split completely in $\frak N$. Now
there are three cases:

\begin{enumerate}
\item
$\left \{ \begin{array}{l}
l \hspace{.05in} \mbox{splits completely in $N_1$,} \\
\mbox{but does not in $N$}
\end{array}
\right \}$ $\iff$ $\left \{ \begin{array}{l}
l \hspace{.05in} \mbox{satisfies} \hspace{.05in} \langle 1, 32\rangle \\
l \hspace{.05in} \mbox{satisfies} \hspace{.05in} \langle p, -2\rangle \\
\end{array}
\right \} $

\item
$\left \{ \begin{array}{l}
l \hspace{.05in} \mbox{splits completely in $N$} \\
\mbox{but does not in $N_1$}
\end{array}
\right \} $ $\iff$ $\left \{ \begin{array}{l}
l \hspace{.05in} \mbox{does not satisfy} \hspace{.05in} \langle 1, 32\rangle \\
l \hspace{.05in} \mbox{satisfies} \hspace{.05in} \langle 1, -2p\rangle \\
\end{array}
\right \} $

\item
$\left \{ \begin{array}{l}
l \hspace{.05in} \mbox{does split completely} \\
\mbox{in $N_1$ or $N$} \\

\end{array}
\right \} $ $\iff$ $\left \{ \begin{array}{l}
l \hspace{.05in} \mbox{does not satisfy} \hspace{.05in} \langle 1, 32\rangle \\
l \hspace{.05in} \mbox{satisfies} \hspace{.05in} \langle p, -2\rangle \\
\end{array}
\right \}.  $

\end{enumerate}
\end{rem}

We can now prove Theorem 1.2

\begin{proof}
Consider the set $ X = \{l \hspace{.05in} \mbox{prime}: l \hspace{.05in}\mbox{is unramified in $\frak N$ and}
\Big(\frac{{\frak N}/\mathbb Q}{l}\Big) = \{k\} \} $ for some $k \in$
Gal$({\frak N}/\mathbb Q)$.  By the $\check C$ebotarev Density Theorem, the set $X$ has natural density
$\frac{1}{32}$ in the set of all primes $l$. Recall $$ A = \{l \hspace{.05in}\mathrm{rational
\hspace{.05in} prime}: l \equiv 1 \bmod 8 \hspace{.05in} \mathrm{and} \hspace{.05in} \Big(\frac{l}{p}\Big) = 1 \} $$ for some fixed prime $p \equiv 1 \bmod 8$.  By Dirichlet's Theorem on primes in arithmetic progressions, $A$ has natural density $\frac{1}{8}$ in the set of
all primes $l$.  Thus $X$ has natural density $\frac{1}{4}$ in $A$. If $p$ satisfies $\langle 1, 32\rangle$, then by Proposition 4.5,
$$
\mbox{4-rank} \hspace{.05in} K_2(\mathcal{O}_{\mathbb Q(\sqrt{pl})}) = 1 \iff
\left \{ \begin{array}{l}
l \hspace{.05in} \mbox{satisfies} \hspace{.05in} \langle 1, 32\rangle \\
l \hspace{.05in} \mbox{satisfies} \hspace{.05in} \langle p, -2\rangle \\
\end{array} \right \} \hspace{.05in} \mbox{or} \hspace{.05in}
\left \{ \begin{array}{l}
l \hspace{.05in} \mbox{does not} \\
\mbox{satisfy} \hspace{.05in} \langle 1, 32\rangle \\
\end{array} \right \}
$$
and
$$
\mbox{4-rank} \hspace{.05in} K_2(\mathcal{O}_{\mathbb Q(\sqrt{pl})}) = 2 \iff
\left \{ \begin{array}{l}
l \hspace{.05in} \mbox{satisfies} \hspace{.05in} \langle 1, 32\rangle \\
l \hspace{.05in} \mbox{satisfies} \hspace{.05in} \langle 1, -2p\rangle \\
\end{array} \right \}.
$$
Using Remark 6.2, we see that for $\mathbb Q(\sqrt{pl})$, 4-rank 1
and 4-rank 2 appear with natural density $\frac{1}{4} +
\frac{1}{2} = \frac{3}{4}$ and $\frac{1}{4}$ respectively. A
similar argument works if $p$ does not satisfy $\langle 1,
32\rangle$. For $\mathbb Q(\sqrt{-pl})$, use Proposition 4.5 and
Remark 6.2 to obtain that 4-rank 1 and 2 each appear with natural
density $\frac{1}{4} + \frac{1}{4} = \frac{1}{2}$ in $A$.
\end{proof}

\section{Proof of Two Corollaries}
For squarefree, odd integers d, recall the sets $X = \{d: d=pl \}$
and $Y = \{d: d=-pl \}$ for distinct primes $p$ and $l$. Now consider the sets
\begin{center}
$X_i = \{d: d=pl, \hspace{.05in} p \equiv i \bmod 8 \}$,

$Y_i = \{d:-pl, \hspace{.05in} p \equiv i \bmod 8 \}.$
\end{center}

Thus $X = X_1 \cup X_3 \cup X_5 \cup X_7$ and $Y = Y_1 \cup Y_3 \cup Y_5 \cup Y_7$. Additionally consider the sets
\begin{center}
$X_{i,j} = \{d: d=pl, \hspace{.05in} p \equiv i \bmod 8, \hspace{.05in} l \equiv j \bmod 8 \}$,

$Y_{i,j} = \{d:d=-pl, \hspace{.05in} p \equiv i \bmod 8, \hspace{.05in} l \equiv j \bmod 8 \}$.
\end{center}

Thus, for example, $X_1 = X_{1,1} \cup X_{1,3} \cup X_{1,5} \cup
X_{1,7}$ and $Y_7 = Y_{7,1} \cup Y_{7,3} \cup Y_{7,5} \cup
Y_{7,7}$.

In Tables 1 and 2 below, for $K_2(\mathcal{O}_{\mathbb
Q(\sqrt{pl})})$, we provide cases in which densities of 4-rank
values follow from congruence conditions on $p$ and $l$, a
condition on the Legendre symbol $\Big(\frac{l}{p}\Big)$ (if any),
and Dirichlet's theorem on primes in arithmetic progressions. In
Tables 3 and 4, we provide the same information for
$K_2(\mathcal{O}_{\mathbb Q(\sqrt{-pl})})$ (compare with
\cite{CHnotes} or Tables I and II in \cite{Qin1} and \cite{Qin2}).

\section*{Table 1: {$\mathbb Q(\sqrt{pl})$}}
\begin{center}
\renewcommand{\arraystretch}{1.75}
\begin{tabular}{|c||c|c|}
\multicolumn{3}{c}{} \\ \hline
$p$, $l$ mod 8 & 4-rank & Densities
\\ \hline
& & \\
3, 3 & 0 &  $\frac{1}{4}$ in $X_3$ \\
& & \\ \hline
& & \\
5, 5 & 1 & $\frac{1}{4}$ in $X_5$ \\
& & \\ \hline
& & \\
7, 7 & 1 & $\frac{1}{4}$ in $X_7$ \\
& & \\ \hline
& & \\
3, 5 & 1 & $\frac{1}{4}$ in $X_3$ and $X_5$ \\
& & \\ \hline
& & \\
3, 7 & 1 & $\frac{1}{4}$ in $X_3$ and $X_7$ \\
& & \\ \hline
& & \\
5, 7 & 1 & $\frac{1}{4}$ in $X_5$ and $X_7$ \\
& & \\ \hline
\end{tabular}
\end{center}

\section*{Table 2: {$\mathbb Q(\sqrt{pl})$}}
\begin{center}
\renewcommand{\arraystretch}{1.75}
\begin{tabular}{|c||c|c|c|}
\multicolumn{4}{c}{} \\ \hline
$p$, $l$ mod 8 & Legendre symbols &
4-rank & Densities \\ \hline
& & & \\
1, 3 & $\big(\frac{l}{p}\big) = -1$ & 0 & $\frac{1}{8}$ in $X_1$ and $X_3$ \\
     & & & \\
     & $\big(\frac{l}{p}\big) = 1$ & 1 & $\frac{1}{8}$ in $X_1$ and $X_3$ \\
& & & \\ \hline
& & & \\
1, 5 & $\big(\frac{l}{p}\big) = -1$ & 0 & $\frac{1}{8}$ in $X_1$ and $X_5$ \\
     & & & \\
     & $\big(\frac{l}{p}\big) = 1$ & 1 & $\frac{1}{8}$ in $X_1$ and $X_5$ \\
& & & \\ \hline
& & & \\
1, 7 & $\big(\frac{l}{p}\big) = -1$ & 1 & $\frac{1}{8}$ in $X_1$ and $X_7$ \\
     & & & \\
     & $\big(\frac{l}{p}\big) = 1$ & 1 & $\frac{1}{16}$ in $X_1$ and $X_7$ \\
     & & 2 & $\frac{1}{16}$ in $X_1$ and $X_7$ \\
& & & \\ \hline
\end{tabular}
\end{center}

\section*{Table 3: {$\mathbb Q(\sqrt{-pl})$}}
\begin{center}
\renewcommand{\arraystretch}{2}
\begin{tabular}{|c||c|c|}
\multicolumn{3}{c}{} \\ \hline
$p$,$l$ mod 8 & 4-rank & Densities
\\ \hline
& & \\
3, 3 & 1 &  $\frac{1}{4}$ in $Y_3$ \\
& & \\ \hline
& & \\
5, 5 & 1 & $\frac{1}{4}$ in $Y_5$ \\
& & \\ \hline
& & \\
7, 7 & 1 & $\frac{1}{4}$ in $Y_7$ \\
& & \\ \hline
& & \\
3, 5 & 0 & $\frac{1}{4}$ in $Y_3$ and $Y_5$ \\
& & \\ \hline
& & \\
3, 7 & 0 & $\frac{1}{4}$ in $Y_3$ and $Y_7$ \\
& & \\ \hline
& & \\
5, 7 & 0 & $\frac{1}{4}$ in $Y_5$ and $Y_7$ \\
& & \\ \hline
\end{tabular}
\end{center}

\section*{Table 4: {$\mathbb Q(\sqrt{-pl})$}}
\begin{center}
\renewcommand{\arraystretch}{2}
\begin{tabular}{|c||c|c|c|}
\multicolumn{4}{c}{} \\
\hline $p$, $l$ mod 8 & Legendre symbols &
4-rank & Densities \\ \hline
& & & \\
1, 1 & $\big(\frac{l}{p}\big) = -1$ & 1 & $\frac{1}{8}$ in $Y_1$. \\
& & & \\ \hline
& & & \\
1, 3 & $\big(\frac{l}{p}\big) = -1$ & 0 & $\frac{1}{8}$ in $Y_1$ and $Y_3$ \\
     & & & \\
     & $\big(\frac{l}{p}\big) = 1$ & 1 & $\frac{1}{8}$ in $Y_1$ and $Y_3$ \\
& & & \\ \hline
& & & \\
1, 5 & $\big(\frac{l}{p}\big) = -1$ & 0 & $\frac{1}{8}$ in $Y_1$ and $Y_5$ \\
     & & & \\
     & $\big(\frac{l}{p}\big) = 1$ & 1 & $\frac{1}{8}$ in $Y_1$ and $Y_5$ \\
& & & \\ \hline
& & & \\
1, 7 & $\big(\frac{l}{p}\big) = -1$ & 0 & $\frac{1}{8}$ in $Y_1$ and $Y_7$ \\
     & & & \\
     & $\big(\frac{l}{p}\big) = 1$ & 0 & $\frac{1}{16}$ in $Y_1$ and $Y_7$ \\
     & & 1 & $\frac{1}{16}$ in $Y_1$ and $Y_7$ \\
& & & \\ \hline
\end{tabular}
\end{center}

\begin{rem} By Theorem 1.2, $p \equiv l \equiv 1 \bmod 8$ with $\Big(\frac{l}{p}\Big) = 1$ yields 4-rank 1 and 2 with densities $\frac{3}{32}$ and $\frac{1}{32}$ respectively in $X_1$. By Theorem 1.3, $p \equiv l \equiv 1 \bmod 8$ with $\Big(\frac{l}{p}\Big) = -1$ yields 4-rank 0 and 1 each with density $\frac{1}{16}$ in $X_1$. We can now prove Corollary 1.4.
\end{rem}

\begin{proof}
Regarding the set $X_1$:
\begin{itemize}
\item 4-rank 0, 1, and 2 appear with natural densities $\frac{1}{16}$, $\frac{3}{32} + \frac{1}{16} = \frac{5}{32}$, and $\frac{1}{32}$ in $X_{1,1}$

\item 4-rank 0 and 1 each appear with natural densities $\frac{1}{8}$ in $X_{1,3}$

\item 4-rank 0 and 1 each appear with natural densities $\frac{1}{8}$ in $X_{1,5}$

\item 4-rank 1 and 2 appear with natural densities $\frac{1}{8} + \frac{1}{16} = \frac{3}{16}$ and $\frac{1}{16}$ in $X_{1,7}$.

\end{itemize}

Thus 4-rank 0, 1, and 2 appear with natural densities
$\frac{5}{16}$, $\frac{19}{32}$, and $\frac{3}{32}$ in $X_1$. For
the set $X_3$:
\begin{itemize}
\item 4-rank 0 and 1 each appear with natural density $\frac{1}{8}$ in $X_{3,1}$
\item 4-rank 0 appears with natural density $\frac{1}{4}$ in $X_{3,3}$
\item 4-rank 1 appears with natural density $\frac{1}{4}$ in $X_{3,5}$
\item 4-rank 1 appears with natural density $\frac{1}{4}$ in $X_{3,7}$.
\end{itemize}
So 4-rank 0 and 1 appear with natural densities $\frac{3}{8}$ and
$\frac{5}{8}$ in $X_3$. Similarly, 4-rank 0 and 1 appear with
natural densities $\frac{1}{8}$ and $\frac{7}{8}$ in $X_5$ and
4-rank 1 and 2 appear with natural densities $\frac{15}{16}$ and
$\frac{1}{16}$ in $X_7$. As each $X_i$ has density $\frac{1}{4}$
in $X$,
\begin{itemize}
\item 4-rank 0 appears with natural density $\frac{5}{64} + \frac{3}{32} + \frac{1}{32} = \frac{13}{64}$ in $X$
\item 4-rank 1 appears with natural density $\frac{19}{128} + \frac{5}{32} + \frac{7}{32} + \frac{15}{64} = \frac{97}{128}$ in $X$
\item 4-rank 2 appears with natural density $\frac{3}{128} + \frac{1}{64} = \frac{5}{128}$ in $X$.
\end{itemize}
\end{proof}

\begin{rem} By Theorem 1.2, $p \equiv l \equiv 1 \bmod 8$ with $\Big(\frac{l}{p}\Big) = 1$ yields 4-rank 1 and 2 each with density $\frac{1}{16}$ in $Y_1$. We can now prove Corollary 1.5.
\end{rem}

\begin{proof}
Regarding the set $Y_1$:
\begin{itemize}
\item 4-rank 1 and 2 appear with natural densities $\frac{1}{8} + \frac{1}{16} = \frac{3}{16}$  and $\frac{1}{16}$ in $Y_{1,1}$

\item 4-rank 0 and 1 each appear with natural densities $\frac{1}{8}$ in $Y_{1,3}$

\item 4-rank 0 and 1 each appear with natural densities $\frac{1}{8}$ in $Y_{1,5}$

\item 4-rank 0 and 1 appear with natural densities $\frac{1}{8}$ and $\frac{1}{16} + \frac{1}{16} = \frac{1}{8}$ in $Y_{1,7}$.
\end{itemize}

Thus 4-rank 0, 1, and 2 appear with natural densities
$\frac{3}{8}$, $\frac{9}{16}$, and $\frac{1}{16}$ in $Y_1$. For
the set $Y_3$:
\begin{itemize}
\item 4-rank 0 and 1 each appear with natural density $\frac{1}{8}$ in $Y_{3,1}$
\item 4-rank 1 appears with natural density $\frac{1}{4}$ in $Y_{3,3}$
\item 4-rank 0 appears with natural density $\frac{1}{4}$ in $Y_{3,5}$
\item 4-rank 0 appears with natural density $\frac{1}{4}$ in $Y_{3,7}$.
\end{itemize}
So 4-rank 0 and 1 appear with natural densities $\frac{5}{8}$ and
$\frac{3}{8}$ in $Y_3$. Similarly, 4-rank 0 and 1 appear with
natural densities $\frac{5}{8}$ and $\frac{3}{8}$ in $Y_5$ and
4-rank 0 and 1 appear with natural densities $\frac{11}{16}$ and
$\frac{5}{16}$ in $Y_7$. As each $Y_i$ has density $\frac{1}{4}$
in $Y$,
\begin{itemize}
\item 4-rank 0 appears with natural density $\frac{3}{32} + \frac{5}{32} + \frac{5}{32} + \frac{11}{64} = \frac{37}{64}$ in $Y$
\item 4-rank 1 appears with natural density $\frac{9}{64} + \frac{3}{32} + \frac{3}{32} + \frac{5}{64} = \frac{13}{32}$ in $Y$
\item 4-rank 2 appears with natural density $\frac{1}{64}$ in $Y$.
\end{itemize}
\end{proof}

\section*{Appendix}
The approach of Hurrelbrink and Kolster in \cite {HK98} led us to
write a program in GP/PARI \cite{pari} which generates the
numerical values in Tables 5-8. The aim is to motivate possible
density results for tame kernels of quadratic number fields. In
Tables 5 and 6, $p$, $l$, and $r$ are distinct odd primes. In
Tables 7 and 8, $d$ is odd and squarefree.

\section*{Table 5}
\begin{center}
\begin{tabular}{|c||c|c|}
\multicolumn{3}{c}{} \\ \hline
Cardinality & $105 \leq d=plr < 10^6$ & $\%$ \\ \hline
$\mid$ 4-rank 0 $\mid$ & 8247 &  6.827 \\ \hline
$\mid$ 4-rank 1 $\mid$ & 92544 &  76.605 \\ \hline
$\mid$ 4-rank 2 $\mid$ & 20000 &  16.555  \\ \hline
$\mid$ 4-rank 3 $\mid$ & 16 &  .013 \\ \hline
\end{tabular}
\end{center}

\section*{Table 6}
\begin{center}
\begin{tabular}{|c||c|c|}
\multicolumn{3}{c}{} \\ \hline
Cardinality & $-10^6 < d=-plr \leq -105$ & $\%$ \\ \hline
$\mid$ 4-rank 0 $\mid$ & 67970 &  56.2633 \\ \hline
$\mid$ 4-rank 1 $\mid$ & 50147 &  41.5100  \\ \hline
$\mid$ 4-rank 2 $\mid$ &  2688 &  2.2250  \\ \hline
$\mid$ 4-rank 3 $\mid$ & 2 & .0017 \\ \hline
\end{tabular}
\end{center}

\section*{Table 7}
\begin{center}
\begin{tabular}{|c||c|c|}
\multicolumn{3}{c}{} \\ \hline
Cardinality & $3 \leq d < 10^6$ & $\%$ \\ \hline
$\mid$ 4-rank 0 $\mid$ & 93736 &  23.1284 \\ \hline
$\mid$ 4-rank 1 $\mid$ & 278138 &  68.6278  \\ \hline
$\mid$ 4-rank 2 $\mid$ &  33148 &  8.1789  \\ \hline
$\mid$ 4-rank 3 $\mid$ & 263 & .0649 \\ \hline
\end{tabular}
\end{center}

\section*{Table 8}
\begin{center}
\begin{tabular}{|c||c|c|}
\multicolumn{3}{c}{} \\ \hline
Cardinality & $-10^6 < d \leq -3$ & $\%$ \\ \hline
$\mid$ 4-rank 0 $\mid$ & 251884 &  62.14985 \\ \hline
$\mid$ 4-rank 1 $\mid$ & 148669 &  36.68258  \\ \hline
$\mid$ 4-rank 2 $\mid$ &  4730 &  1.16708  \\ \hline
$\mid$ 4-rank 3 $\mid$ & 2 & .00049 \\ \hline
\end{tabular}
\end{center}

\section*{acknowledgments} We thank the referee for informing us of the paper by Qin Yue.  We also thank Manfred Kolster for his suggestions and comments.

\end{document}